\documentclass{amsart}
\usepackage{amssymb, amsthm}
\usepackage{epsfig}

\def\ovl{\overline}
\def\ol{\overline}
\def\N{\mathbb{N}}
\def\Z{\mathbb{Z}}
\def\R{\mathbb{R}}
\def\C{\mathbb{C}}
\def\disk{\mathbb D}
\def\diskbar{\ovl{\disk}}
\def\Cc{\widehat{{\C}}}
\def\H{\mathbb{H}}

\def\Re{{\rm Re}}
\def\Im{\mbox{\rm Im}}
\def\id{\mbox{\rm id}}
\def\tilde{\widetilde}
\def\area{\mbox{area}}
\def\mod{\mbox{\rm mod}}

\def\sm{\setminus}
\def\phi{\varphi}
\def\eps{\varepsilon}

\newtheorem{theorem}{Theorem}[section]
\newtheorem{proposition}[theorem]{Proposition} 
\newtheorem{corollary}[theorem]{Corollary}
\newtheorem{lemma}[theorem]{Lemma}

\theoremstyle{definition}
\newtheorem{definition}[theorem]{Definition}
\theoremstyle{remark}
\newtheorem{remark}[theorem]{Remark}

\title{On Newton's Method for Entire Functions}

\author{Johannes R\"uckert}
\address{School of Engineering and Science, International University Bremen, Postfach 750 561,  28725 Bremen, Germany}
\email{j.rueckert@iu-bremen.de}

\author{Dierk Schleicher}
\address{School of Engineering and Science, International University Bremen, Postfach 750 561,  28725 Bremen, Germany}
\email{dierk@iu-bremen.de}

\begin{document}
\maketitle
\begin{abstract}
    \noindent
    The Newton map $N_f$ of an entire function $f$ turns the roots of $f$
    into attracting fixed points. Let $U$ be the immediate
    attracting basin for such a fixed point of $N_f$.

    We study the behavior of $N_f$ in a component $V$ of $\C\sm
    U$. If $V$
    can be surrounded by an invariant curve within $U$ and
    satisfies the condition that for all $z\in\Cc$, $N_f^{-1}(\{z\})\cap V$ is a finite
    set, we show that $V$ contains another immediate basin of $N_f$ or a {\em virtual immediate
    basin} (Definition \ref{Def_VirtualBasins}).
\end{abstract}

\section{Introduction}
\noindent Newton's method is a classical way to approximate roots
of entire functions by an iterative procedure. Trying to
understand this method may very well be called the founding
problem of holomorphic dynamics \cite[p.~51]{Milnor}.

Newton's method for a complex polynomial $p$ is the iteration of a
rational function $N_p$ on the Riemann sphere. Such dynamical
systems have been extensively studied in recent years. Tan Lei
\cite{Tan} gave a complete classification of Newton maps of cubic
polynomials. In 1992, Manning \cite{Manning} constructed a finite
set of starting values for $N_p$ that depends only on the degree
of $p$, such that for any appropriately normalized polynomial with
degree $d\geq 10$, the set contains at least one point that
converges to a root of $p$ under iteration of $N_p$. Hubbard,
Schleicher and Sutherland \cite{HSS} extended this by constructing
a small set of starting values that depends only on the degree
$d\geq 2$ and trivial normalizations and finds all roots of $p$.

\begin{figure}[htb]
    \begin{center}
        \setlength{\unitlength}{1cm}
        \begin{picture}(6,6)
            \epsfig{file=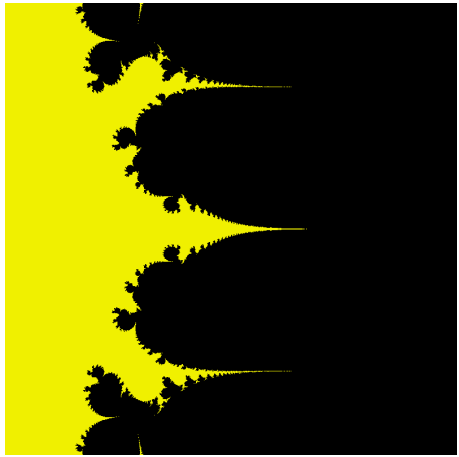, width=6cm}
        \end{picture}
        \caption{Newton's map for $z\mapsto z\,e^{e^z}$. The immediate basin of $0$
        has infinitely many accesses to
        the right. Any two of them surround a virtual immediate basin.
             More precisely, all curves of the form $(2k+1)\pi i + [2,\infty]$ are contained in a virtual immediate basin; the virtual basins for $k_1\neq k_2$ are disjoint and separated by an access to $\infty$ of the immediate basin of $0$.     
The visible area is from $-8-10i$ to $12+10i$.}
        \label{Figure_SebastianDiplom}
    \end{center}
\end{figure}

If $f$ is a transcendental entire function, the associated Newton
map $N_f$ will generally be transcendental meromorphic, except in
the special  case $f=pe^q$ with polynomials $p$ and $q$ (see
Proposition \ref{Prop_RationalNewton}) which was studied by Haruta
\cite{Haruta}. Bergweiler \cite{Bergweiler2} proved a
no-wandering-domains theorem for transcendental Newton maps that
satisfy several finiteness assumptions. Mayer and Schleicher
\cite{MS} have shown that immediate basins for Newton maps of
entire functions are simply connected and unbounded, extending a
result of Przytycki \cite{Przytycki} in the polynomial case. They
have also shown that Newton maps of transcendental functions may
exhibit a type of Fatou component that does not appear for Newton
maps of polynomials, so called {\em virtual immediate basins}
(Definition \ref{Def_VirtualBasins}) in which the dynamics
converges to $\infty$. The thesis \cite{SebastianDiplom}
investigates the Newton map of the transcendental function
$z\mapsto z e^{e^z}$ and shows that it exhibits virtual immediate
basins; see Figure \ref{Figure_SebastianDiplom} for an
illustration. While immediate basins of roots are by definition
related to zeroes of $f$ (compare Definition
\ref{Def_ImmediateBasin}), under mild technical assumptions a virtual immediate basin leads to an {\em asymptotic} zero of $f$; in other words, a virtual immediate basin often contains
an asymptotic path of an asymptotic value at $0$ for $f$ \cite{Buff}.

In this paper, we continue the work of \cite{MS} and investigate
the behavior of Newton maps in the complement of an immediate
basin. Our main result (Theorem \ref{Thm_NecessaryCondition}) is
that if a complementary component can be surrounded by an
invariant curve through $\infty$, then it contains another
immediate basin or virtual immediate basin, unless it maps infinite-to-one
onto at least one point of $\Cc$. We believe that the last ``unless''-condition is
unnecessary, but our methods do not allow us to show this.

An immediate corollary for
Newton maps of polynomials is that between any two ``channels'' of any root, there
is always another root. This is folklore, but we do not know of a published
reference. This result can be viewed as a first step
towards a classification of polynomial Newton maps.

Our paper is structured in the following way: In Section
\ref{Sec_Newton}, we give an introduction to some general
properties of Newton maps. In Section \ref{Sec_Curves}, we
investigate homotopy classes of curves to $\infty$ in immediate
basins and prove some auxiliary results. In Section
\ref{Sec_Lefschetz}, we prove a fixed point theorem which we will
need and which might be interesting in its own right. In Section
\ref{Sec_Main}, we state and prove our main result.

\section{Newton's Method as a Dynamical System}
\label{Sec_Newton}

\subsection{Immediate Basins}
Let $f:\C\to\C$ be a non-constant entire function and $N_f$ its
associated (meromorphic) Newton map
\[
    N_f \;=\; \id\,-\,\frac{f}{f'}\;\;.
\]
If $f$ is a polynomial, then $N_f$ extends to a rational map
$\Cc\to\Cc$. If $\xi$ is a root of $f$ with multiplicity $m\geq
1$, then it is an attracting fixed point of $N_f$ with multiplier
$\frac{m-1}{m}$. Conversely, every fixed point $\xi\in\C$ of $N_f$
is attracting and a root of $f$.
\begin{definition}[Immediate Basin]
\label{Def_ImmediateBasin}
    Let $\xi$ be an attracting fixed point of $N_f$. The {\em basin of $\xi$}
    is $\{z\in\C:\,\lim_{n\rightarrow\infty} N_f^{\circ n}(z) =
    \xi\}$, the open set of points
    which converge to $\xi$ under iteration. The
    connected component $U$ of the basin that contains $\xi$ is called its \emph{immediate
    basin}.
\end{definition}
Immediate basins are $N_f$-invariant because they are Fatou
components and contain a fixed point. The following theorem is the
main result (Theorem 2.7) of \cite{MS}.
\begin{theorem}[Immediate Basins Simply Connected]
\label{Thm_ImmediateBasinSC}
    If $\xi$ is an attracting fixed point of the Newton map $N_f$, then its
    immediate basin $U$ is simply connected and
    unbounded.
    \qed
\end{theorem}
We will use the following notation throughout the paper:

If $\gamma$ is a curve, the symbol $\gamma$ denotes the mapping
$\gamma:I\to\C$ from an interval into the plane as well as its
image $\gamma(I)\subset\C$. By a {\em tail} of an unbounded curve
we mean any unbounded connected part of its image.

For $r>0$ and $z\in\C$, the symbol $B_r(z)$ designates the disk of
radius $r$ centered at $z$.

The {\em full preimage} of a point $z\in\Cc$ is the set
$N_f^{-1}(\{z\})$. Its only accumulation point can be $\infty$ by
the identity theorem. Any point $z'\in N_f^{-1}(\{z\})$ is called
a {\em preimage} of $z$.

Unless stated otherwise, the boundary and the closure of a set are
considered in $\C$.

\subsection{Singular Values}
Since the concept of singular values is crucial for the study of
dynamical systems, we give a brief reminder of the most important
types. In particular, we state some properties of {\em asymptotic
values}; these appear only for transcendental maps.

\begin{definition}[Singular Value]
Let $h:\C\to\Cc$ be a meromorphic function.  We call a point $p\in
\C$ a {\em regular point} of $h$ if $p$ has a neighborhood on
which $h$ is injective. Otherwise, we call $p$ a {\em critical
point}. A point $v\in\Cc$ is called a {\em regular value} if there
exists a neighborhood $V$ of $v$ such that for every component $W$
of $h^{-1}(V)$, $h^{-1}|_V:V\to W$ is a single-valued meromorphic
function. Otherwise, $v$ is called a {\em singular value}.

The image of a critical point is a singular value and is called a
{\em critical value}.
\end{definition}
Critical points in $\C$ are exactly the zeroes of the first
derivative. For a rational map, all singular values are critical
values.
\begin{definition}[Asymptotic Value]
    Let $h:\C\to\Cc$ be a transcendental meromorphic function. A point $a\in\Cc$ is
    called an \emph{asymptotic value} of $h$ if there exists a curve
    $\Gamma:\R_+\to\C$ with $\lim_{t\to\infty}\Gamma(t)=\infty$ such that
    $\lim_{t\rightarrow\infty}h(\Gamma(t))=a$.
    We call $\Gamma$ an \emph{asymptotic path} of $a$.
\end{definition}
In general, an asymptotic value is defined by having an asymptotic
path towards any essential singularity. Note that in our
definition, the set of singular values is the closure of the set
of critical and asymptotic values.

We follow \cite{BergweilerEremenko} in the classification of
asymptotic values.
\begin{definition}[Direct and Indirect Singularity]
    Let $h:\C\to\Cc$ be a meromorphic function and $a\in\C$ be a
    finite asymptotic value with asymptotic path $\Gamma$.
    For each $r>0$, let $U_r$ be the unbounded component of
    $h^{-1}(B_r(a))$ that contains an unbounded end of $\Gamma$.

    We say that $a$ is a \emph{direct singularity} (with respect
    to $\Gamma$) if there is an $r>0$ such that $h(z)\neq
    a$ for all $z\in U_r$. We call $a$ an \emph{indirect singularity}
    if for all $r>0$, there is a $z\in U_r$ such that $h(z)=a$
    (then there are infinitely many such $z$ in $U_r$).
\end{definition}
\begin{theorem}[Direct Singularities]
\label{Thm_DirectSingularities} \emph{\cite[Theorem 5]{Heins}}.
The set of direct singularities of a meromorphic function is
always countable.
\qed
\end{theorem}
It is possible however that the set of (direct and indirect)
singularities is the entire extended plane:
Eremenko~\cite{Eremenko} constructed meromorphic functions of prescribed
finite order whose set of asymptotic values is all of $\Cc$.
\begin{lemma}[Unbounded Preimage]
\label{Lem_AsymptoticValueOnBoundary}
    Let $h:\C\to\Cc$ be a meromorphic function and $B\subset\C$ a bounded topological
    disk whose boundary is a simple closed curve $\beta$. Suppose that $\beta$ contains no critical
    values and that $\tilde{B}$ is an unbounded preimage component of
    $B$. Then $\partial\tilde{B}$ contains an unbounded curve
    $\tilde{\beta}$ with $h(\tilde{\beta})\subset\beta$ such that
    either $h|_{\tilde{\beta}}:\tilde{\beta}\to \beta$ is a universal
    covering map or $h(\tilde{\beta})$ lands at an asymptotic
    value on $\beta$.
\end{lemma}
\begin{proof} Let $w\in\partial\tilde{B}$. Clearly, $h(w)\in\beta$ and by
assumption, $h$ is a local homeomorphism  in a neighborhood of
$w$. It follows that the closed and unbounded set
$\partial\tilde{B}$ is locally an arc everywhere; therefore it
cannot accumulate in any compact subset of $\C$ and must contain
an arc $\tilde{\beta}$ that converges to $\infty$. The curve
$\tilde{\beta}$ contains no critical points. If
$h|_{\tilde{\beta}}:\tilde{\beta}\to \beta$ is not a universal
covering map, then it must land at an asymptotic value.\end{proof}

\subsection{Newton Maps}
We show that there is only one class of entire functions that have
rational Newton maps. This class contains all polynomials. We give
a classification of the dynamics within immediate basins for
Newton maps of polynomials.

First, we investigate under which conditions a meromorphic
function is the Newton map of an entire function. The following
proposition uses ideas of Matthias G\"orner and extends a similar
result for rational maps (see below) and certain transcendental
functions \cite[page 3]{Bergweiler2}. We do not know if the proposition is new; however, we certainly do not know of a published reference.
\begin{proposition}[Newton Maps]
\label{Prop_NewtonMaps} Let $N:\C\to\Cc$ be a meromorphic
function. It is the Newton map of an entire function $f:\C\to\C$
if and only if for each fixed point $N(\xi)=\xi\in\C$, there is a
natural number $m\in\N$ such that $N'(\xi)=\frac{m-1}{m}$. In this
case, there exists $c\in\C\sm\{0\}$ such that
\[
    f = c\cdot \exp\left({\int\frac{d\zeta}{\zeta-N(\zeta)}}\right)\;\;.
\]
Two entire functions $f,g$ have the same Newton maps if and only
if $f=c\cdot g$ for a constant $c\in\C\sm\{0\}$.
\end{proposition}
\begin{proof} We start with the last claim: $f$ and $cf$ have the same
Newton map $\id-f/f'=\id-1/(\ln f)'$. Conversely, if $f$ and $g$
have the same Newton maps, then $(\ln f)'=(\ln g)'$, and the claim
follows.

It is easy to check that every Newton map satisfies the
criterion on derivatives at fixed points.

For the other direction, we construct a map $f$ such that
$N_f = N$. Let $z_0\in\C$ be any base point and define
$\tilde{f}(z)=\int_\gamma \frac{d \zeta}{\zeta-N(\zeta)}$, where
$\gamma:[0;1]\to\C$ is any integration path from $z_0$ to $z$ that
avoids the fixed points of $N$. This defines $\tilde{f}$ up to
$2\pi i k$: if $\gamma'$ is another choice of integration path,
the residue theorem shows that
\[
    \frac{1}{2 \pi i}\int_{\gamma'\circ\gamma^{-1}} \frac{d \zeta}{\zeta-N(\zeta)}
= \sum_{N(\xi)=\xi}
    \mbox{Res}_{\xi} \left(\frac{1}{\zeta-N(\zeta)}\right)\;\;,
\]
where the sum is taken over the finitely many fixed points of $N$
that are contained in the compact regions bounded by the closed
path $\gamma'\circ\gamma^{-1}$. Near a fixed point $\xi$, it is
easy to show that $z-N(z)=\frac{1}{m}(z-\xi)+o(z-\xi)$. Hence we
get $\mbox{Res}_{\xi}\left(\frac{1}{z-N(z)}\right)=m\in\N$.

It follows that the map $f=\exp(\tilde{f})$ is well defined and
holomorphic outside the fixed points of $N$. Near such a fixed
point $\xi$, $\tilde{f}$ has the form $m\log(z-\xi)+O(1)$.
Clearly, the real part of this converges to $-\infty$ for
$z\to\xi$, hence setting $f(\xi)=0$ makes $f$ an entire function
as desired. An easy calculation then shows that $N_f = N$. A
different choice of base point $z_0$ will change $f$ by a
multiplicative constant and lead to the same Newton map $N_f$.
\end{proof}
The following corollary is essentially due to Janet Head
(\cite[Proposition 2.1.2]{head}, \cite[Lemma 2.2]{Tan}).
\begin{corollary}[Rational Newton Maps]
\label{HeadTheorem} A rational map $f:\Cc\rightarrow\Cc$ of degree
$d\geq 2$ is the Newton map of a
    polynomial of degree at least two if and only if $f(\infty)=\infty$ and for all other
    fixed points $a_1,\ldots,a_d\in\C$ there exists a number
    $m_j\in\N$ such that $f'(a_j)=\frac{m_j-1}{m_j}<1$.
    Then, $f$ is the Newton map of the polynomial
    \[
        p(z)=a\prod_{j=1}^d (z-a_j)^{m_j}
    \]
    for any complex $a\neq 0$.
\end{corollary}
\begin{proof}[Sketch of Proof: ] Let $a\in\C\sm\{0\}$. Since $N_p$ and $f$ have the same
fixed points with identical multiplicities, the residuals of the
maps $\tilde{f}:=(f-\id)^{-1}$ and $\tilde{N}:=(N_p-\id)^{-1}$ at
their common simple poles $a_1,\dots,a_d\in\C$ agree, and thus
also those at $\infty$. Hence, $\tilde{f}-\tilde{N}$ is a
polynomial with $\lim_{z\to\infty}(\tilde{f}-\tilde{N})(z)=0$.
Hence $\tilde{f}=\tilde{N}$ and the claim follows.
\end{proof}
We want to exclude the trivial case of Newton maps with degree
one.
\begin{lemma}[One Root]
Let $f:\C\to\C$ be an entire function such that its Newton map
$N_f$ has an attracting fixed point $\xi\in\C$ with immediate
basin $U=\C$. Then, there exist $d>0$ and $a\in\C$ such that
$f(z)=a(z-\xi)^d$.
\end{lemma}
\begin{proof} Since $N_f$ has no periodic points of minimal period at
least $2$, it cannot be transcendental \cite[Theorem
2]{Bergweiler}. Hence $N_f$ is rational and its fixed points can
only be $\xi$ and $\infty$, both of which must be simple. It
follows that $N_f$ has degree at most one and since it has no
poles in $\C$, it is a polynomial. The claim now follows from
Proposition \ref{Prop_NewtonMaps}.
\end{proof}

In the rest of this paper, we will assume that $N_f$ is not a
M\"obius transformation. Theorem \ref{Thm_ImmediateBasinSC}
implies then that for each immediate basin $U$ of $N_f$, there
exists a Riemann map $\phi:\disk\to U$ with $\phi(0)=\xi$.

The following simple proposition classifies rational Newton maps
of entire functions. Its first half is stated without proof in
\cite{Bergweiler}.
\begin{proposition}[Rational Newton Map]
\label{Prop_RationalNewton}
    Let $f:\C\to\C$ be an entire function. Its Newton
    map $N_f$ is rational if and only if there are polynomials $p$ and $q$ such that $f$ has the form
    $f=p\,e^q$. In this case, $\infty$ is a repelling or parabolic
    fixed point.

    More precisely, let $m,n\geq 0$ be the degrees of $p$ and $q$,
    respectively. If $n=0$ and $m\geq 2$, then $\infty$ is
    repelling with multiplier $\frac{m}{m-1}$. If $n=0$ and
    $m=1$, then $N_f$ is constant. If $n>0$, then $\infty$ is
    parabolic with multiplier $+1$ and multiplicity $n+1\geq 2$.
\end{proposition}
\begin{proof} By \cite[Corollary 12.7]{Milnor}, every rational function of
degree at least $2$ has a repelling or parabolic fixed point.
Since $N_f$ is a Newton map, this non-attracting fixed point is
unique and must be at $\infty$. In addition to this, there are
finitely many attracting fixed points $a_1,\ldots,a_n\in\C$ with
associated natural numbers $m_1,\ldots,m_n\in\N$ such that the
multipliers satisfy $N_f'(a_i)=\frac{m_i-1}{m_i}$. Let
$p(z)=\prod_{i=1}^n (z-a_i)^{m_i}$.

Since attracting fixed points of $N_f$ correspond exactly to the
roots of $f$, $f$ has the form $f=p\,e^h$ for an entire function
$h$. If $h$ was transcendental, so would be
\[
    N_f\,=\,\id -\frac{p\,e^h}{p'e^h+h'p\,e^h}\;=\;\id-\frac{p}{p'+h'p}\;\;,
\]
a contradiction. The other direction follows by direct calculation
and the rest of the proof is left to the reader.
\end{proof}

\begin{figure}[hbt]
    \begin{center}
        \setlength{\unitlength}{1cm}
        \begin{picture}(6,6)
            \epsfig{file=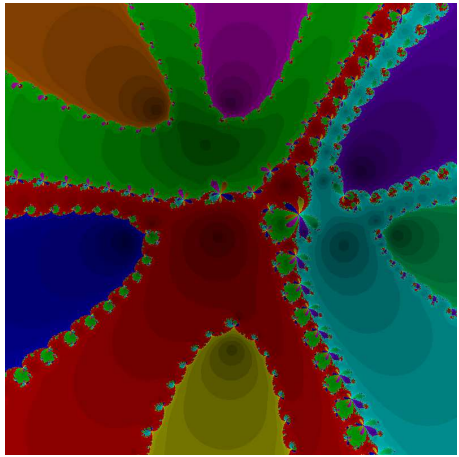, width=6cm}
        \end{picture}
        \caption{Newton map for a polynomial of degree $9$. The channels are clearly visible.}
        \label{Figure_PolyNewton}
    \end{center}
\end{figure}

For Newton maps of $f=pe^q$, the area of every immediate basin is finite if $\deg
q\ge 3$ \cite{Haruta} and infinite if $p(z)=z$ and $\deg q\in\{0,1\}$ \cite{Figen}.

The dynamics within immediate basins of Newton maps of polynomials
has an easy classification, because all singular values are
critical values.
\begin{theorem}[Polynomial Newton Maps]
\label{Thm_PolyNewton} {\em \cite{HSS}} Let $p$ be a polynomial of
degree $d>1$, normalized so that its roots are contained in the
unit disk $\disk$. Let $\xi$ be a root of $p$ and $U$ its
immediate basin for $N_p$. Then, $U$ contains $0<k<d$ critical
points of $N_p$ and $N_p|_U$ is a proper self-map of degree $k+1$.
Outside the disk $B_2(0)$, $N_f$ is conformally conjugate to
multiplication by $\frac{d-1}{d}$. Finally, $U\sm B_2(0)$ has
exactly $k$ unbounded components, so called {\em channels}, each
of which maps over itself under $N_f$.
\qed
\end{theorem}
Figure \ref{Figure_PolyNewton} illustrates this theorem.

\section{Accesses in Immediate Basins}
\label{Sec_Curves}
\subsection{Invariant Accesses}
We investigate the immediate basins of attraction for the
attracting fixed points of $N_f$. If $f$ is a polynomial, we have
seen in Theorem \ref{Thm_PolyNewton} that immediate basins have an
easy geometric structure. In the general case, $N_f$ has an
essential singularity at $\infty$ and immediate basins may well
have infinitely many accesses to $\infty$. We use prime end theory
to distinguish them.

Under a finiteness assumption, we have some control over the image
of a sequence that converges to $\infty$ through an immediate
basin.
\begin{lemma}[Invariant Boundary]
\label{No_Asymptotic_Path}
    Let $U$ be an immediate basin of the Newton map $N_f$ and
    $U_R$ an unbounded component of $U\sm B_R(0)$ with the property that
    no point has infinitely many preimages in $U_R$. Then for any
    sequence $(z_n)\subset U_R$ with $z_n\to\infty$, all limit points of $N_f(z_n)$
    are contained in $\partial U\cup\{\infty\}$.
\end{lemma}
The condition is necessary, because if there exists a point $p\in
U$ with infinitely many preimages $p_1,p_2,\ldots \in U_R$, we
have $p_n\to\infty$ and $N_f(p_n)=p\in U$ for all $n\in\N$.
\begin{proof} Assume there exists a sequence $(z_n)\subset U_R$ that
converges to $\infty$ with $N_f(z_n)\to p\in U$. Let $B\subset U$
be a closed neighborhood of $p$ inside $U$ such that its boundary
$\partial B$ is a simple closed curve $\beta$ that contains no
direct singularities (this is possible by Theorem
\ref{Thm_DirectSingularities}) nor critical values.

Suppose first that $p\not\in N_f(\partial B_R(0))$. Then we may
choose $B$ small enough such that $B\cap N_f(\partial
B_R(0))=\emptyset$. The image of the first finitely many $z_n$
need not be in $B$; ignoring those, each $z_n$ is contained in a
component $W_n$ of $N_f^{-1}(B)\cap U$. If a $W_n$ is bounded, it
maps surjectively onto $B$ under $N_f$. Therefore, by the
finiteness assumption, there can be only finitely many bounded
$W_n$. Each bounded $W_n$ contains finitely many $z_n$; hence
there must be an $n$ such that $W_n$ is unbounded. By Lemma
\ref{Lem_AsymptoticValueOnBoundary} and again because of the
finiteness assumption, $\partial W_n$ contains an asymptotic path
of an asymptotic value on $\beta$. But this asymptotic value must
be an indirect singularity, which also contradicts the finiteness
assumption.

If $p\in N_f(\partial B_R(0))$, a small homotopy of the curve
$\partial B_R(0)$ in a neighborhood of $p$ solves the problem.
\end{proof}
Figure \ref{Figure_SebastianDiplom} suggests that immediate basins
can reach out to infinity in several different directions. We make
this precise in the following definitions that generalize the
concept of a channel in the polynomial case.
\begin{definition}[Invariant Access]
\label{DefFixedAccess}
    Let $\xi$ be an attracting fixed point of $N_f$ and $U$ its immediate
    basin. An \emph{access to $\infty$} of $U$ is a homotopy class of curves within $U$ that begin at
    $\xi$, land at $\infty$ and are homotopic with fixed endpoints.

    An \emph{invariant access to $\infty$} is an access with the additional property that for each representative
    $\gamma$, its image $N_f(\gamma)$ belongs to the access as well.
\end{definition}
\begin{lemma}[Access Induces Prime End]
\label{LemAccessPrimeEnd}
    Let $[\gamma]$ be an access to $\infty$ in $U$. Then
    $[\gamma]$ induces a prime end $\mathcal{P}$ in $U$ with impression
    $\{\infty\}$. If $[\gamma]$ is invariant, then $N_f(\mathcal{P})=\mathcal{P}$.
\end{lemma}
\begin{proof} Let $\gamma\subset U$ be a curve representing $[\gamma]$
that starts at the fixed point $\xi$ and lands at $\infty$. For
$n\in\N$, let $W_n$ be the component of $U\sm B_n(0)$ that
contains a tail of $\gamma$. The $W_n$ represent a prime end
$\mathcal{P}$ with impression $\infty$. Now a curve
$\gamma'\subset U$ that starts at $\xi$ and lands at $\infty$ is
homotopic to $\gamma$ if and only if a tail of it is contained in
$W_n$ for $n$ large enough. Hence the prime end $\mathcal{P}$ of
$[\gamma]$ is well-defined. The last claim follows immediately
from the definition.
\end{proof}

It is clear that different accesses induce different prime ends.
We state one more well-known topological fact about the boundary
behavior of Riemann maps before using prime ends to characterize
invariant accesses.
\begin{lemma}[Accesses Separate Disk]
\label{Lem_RiemannConnected}
    Let $U\subsetneq \C$ be a simply connected unbounded domain
    and $\gamma_1,\gamma_2:\R_0^+\to U\cup\{\infty\}$ two non-homotopic curves that land at $\infty$
    and are disjoint except for their common base point $z_0=\gamma_1(0)=\gamma_2(0)\in U$. Let $C$ be
    a component of $\C\sm (\gamma_1\cup\gamma_2)$ and
    $\phi:\disk\to U$ a Riemann map with $\phi(0)=z_0$.

    Then $\phi^{-1}(\gamma_1)$ and $\phi^{-1}(\gamma_2)$ land at
    distinct points $\zeta_1$ and $\zeta_2$ of $\partial\disk$.
    Furthermore, $\ol{\partial U \cap C}\subset\Cc$ corresponds
    under $\phi^{-1}$ to a closed interval on $\partial\disk$ that
    is bounded by $\zeta_1$ and $\zeta_2$.
\end{lemma}
This follows immediately because $\phi$ extends to a homeomorphism
from $\diskbar$ to the Carath\'eodory compactification of $U$, see
\cite[Theorem 17.12]{Milnor}.

If $f$ is a polynomial, it follows from Theorem
\ref{Thm_PolyNewton} that every immediate basin contains a curve
that lands at $\infty$, is homotopic to its image and induces an
invariant access. In the general case, it is a priori not even
clear that a curve that lands at $\infty$ {\em and} is homotopic
within $U$ to its image induces an invariant access. The following
proposition deals with this issue.
\begin{proposition}[Curve Induces Invariant Access]
\label{PropInvAccesshasFiniteDegree}
    Let $\gamma\subset U\cup\{\infty\}$ be a curve connecting the fixed point $\xi$ to
    $\infty$ such that $N_f(\gamma)$ is homotopic to $\gamma$ in $U$ with
    endpoints fixed. Let $W_n$ be a sequence of fundamental neighborhoods representing
    the prime end $\mathcal{P}$ induced by $[\gamma]$.
    Then $\gamma$ defines an invariant access to
    $\infty$ if and only if there is no $z\in\Cc$ that has infinitely many preimages in all $W_n$.
\end{proposition}
\begin{proof} Suppose that $\gamma$ defines an invariant access, i.e.\ if
$\gamma'$ is homotopic in $U$ to $\gamma$, then $N_f(\gamma')$ is
homotopic to $N_f(\gamma)$. Assume there is a point $z_0\in\Cc$
with the property that $N_f^{-1}(\{z_0\})\cap W_n$ is an infinite
set for all $W_n$. Without loss of generality, we may assume that
for all $n\in\N$, $W_n\sm W_{n+1}$ contains one preimage of
$z_0$. Then we can find a curve $\gamma'$ with a tail contained in
each $W_n$ that goes through a preimage of $z_0$ in each
$W_n\sm W_{n+1}$. Clearly, $\gamma'$ is homotopic to
$\gamma$, while its image does not land at $\infty$ and can
therefore not be homotopic to $N_f(\gamma)$ with endpoints fixed,
a contradiction.

Now suppose that no point has infinitely many preimages in all
$W_n$. Since the $W_n$ are nested, no point can have infinitely
many preimages in any $W_n$ for $n$ sufficiently large. We
uniformize $U$ to the unit disk via a Riemann map $\phi:\disk\to
U$ such that $\phi(0)=\xi$ and consider the induced dynamics
$g=\phi^{-1}\circ N_f\circ \phi:\disk\to \disk$.

By \cite[Corollary 17.10]{Milnor}, $\phi^{-1}(\gamma)$ and
$\phi^{-1}(N_f(\gamma))$ land on $\partial\disk$. Since the curves
are homotopic, they even land at the same point
$\zeta\in\partial\disk$. Now by assumption, there exists an
$\eps>0$ such that within $B_\eps(\zeta)$, no $g$-preimage of any
point in $\disk$ accumulates. By Lemma \ref{No_Asymptotic_Path} it
follows that the $g$-image of any sequence converging to $\partial
\disk$ inside $B_\eps(\zeta)\cap\disk$ will also converge to
$\partial \disk$. Hence we can use the Schwarz Reflection
Principle \cite[Theorem 11.14]{Rudin} to extend $g$
holomorphically to a neighborhood of $\zeta$ in $\C$. It follows
that for the extended map, $\zeta$ is a repelling fixed point with
positive real multiplier: if the multiplier was not positive real,
$g$ would map points in $B_{\eps}(\zeta)\cap\disk$ out of $\disk$.
Also, $\zeta$ cannot be attracting or parabolic, because in this
case it would attract points in $\disk$, which all converge to $0$
under iteration.

Since $\disk$ is simply connected, all curves in $\disk$ from $0$
to $\zeta$ will be homotopic to each other and their $g$-images. A
curve in $\disk$ that starts at $0$ lands at $\zeta$ if and only
if its $\phi$-image in $U$ is homotopic to $\gamma$ with endpoints
fixed, because $\phi^{-1}(\mathcal{P})$ is a prime end in $\disk$
with impression $\zeta$.
\end{proof}
\begin{remark}
We have shown that each invariant access defines a boundary fixed
point in the conjugated dynamics on the unit disk, and the
dynamics can be extended to a neighborhood of this boundary fixed
point, necessarily yielding a repelling fixed point. By
\cite[Corollary 17.10]{Milnor} it follows that different invariant
accesses induce distinct boundary fixed points.

If $f$ is a polynomial, there exists a one-to-one correspondence
between accesses to $\infty$ of $U$ and boundary fixed points of
the induced map $g$ \cite[Proposition 6]{HSS}.
\end{remark}
\begin{corollary}[Invariant Curve]
\label{invariantCurve}
    Each invariant access has an invariant representative,
    i.e.\ a curve $\gamma:\R_0^+\to U$ that lands at $\infty$ with
    $\gamma(0)=\xi$ and $N_f(\gamma)=\gamma$.
\end{corollary}
\begin{proof} For the extension of $g$ to a neighborhood of
$\zeta\in\partial\disk$, the multiplier of $\zeta$ is positive
real. A short piece of straight line in linearizing coordinates
around $\zeta$ maps over itself under $g$. Its forward orbit lands
at the fixed point.
\end{proof}
Since there are uncountably many choices of such invariant curves,
we can always find one that contains no critical or direct
asymptotic values outside a sufficiently large disk.

\subsection{Virtual Basins}
If $f$ is a polynomial and $U\subset\C$ an invariant Fatou
component of $N_f$, then $U$ is the immediate basin of a root of $f$,
because the Julia set of $N_f$ is connected \cite{Shishikura}, all
finite fixed points are attracting and the fixed point at $\infty$
is repelling. If $f$ is transcendental entire, $N_f$ may possess
invariant unbounded Fatou domains in which the dynamics converges
to $\infty$. Such components are  Baker domains or attracting
petals of an indifferent fixed point at infinity. In many cases, such components contain an asymptotic path of an
asymptotic value at 0 for $f$ \cite{Buff}.
\begin{definition}[Virtual Basin]
\label{Def_VirtualBasins}
    An unbounded domain $V\subset\C$ is called \emph{virtual
    immediate basin of $N_f$} if it is maximal (among domains in $\C$) with respect to the
    following properties:
    \begin{enumerate}
        \item $\lim\limits_{n\rightarrow\infty} N_f^{\circ n}(z)=\infty$ for all $z\in
        V$;
        \item \label{property_2} there is a connected and simply connected subdomain $S_0\subset
        V$ such that $N_f(\ol{S_0})\subset S_0$ and for all
        $z\in V$ there is an $m\in\N$ such that $N_f^{\circ
        m}(z)\in S_0$.
    \end{enumerate}
    We call the domain $S_0$ an \emph{absorbing set} for $V$.
\end{definition}
Clearly, virtual immediate basins are forward invariant.
\begin{theorem}[Virtual Basin Simply Connected]
    \emph{\cite[Theorem 3.4]{MS}} Virtual immediate basins are simply
    connected.
    \qed
\end{theorem}
It might be possible to extend Shishikura's theorem
\cite{Shishikura} to show that for Newton maps of entire
functions, all Fatou components are simply connected. Taixes has announced partial results in this direction, in particular he rules out the existence of cycles of Herman rings (see also Corollary \ref{Cor_Herman} below). If it were also known that Baker domains are always simply connected, then a result of Cowen \cite[Theorem 3.2]{Cowen} would imply that every invariant Fatou
component of a Newton map is an immediate basin or a virtual
immediate basin (see \cite[Remark 3.5]{MS}).

\section{A Fixed Point Theorem}
\label{Sec_Lefschetz} Let $X$ be a compact, connected and
triangulable real $n$-manifold and let $f:X\to X$ be continuous
with finitely many fixed points.
Each fixed point of $f$ has a well-defined {\em Lefschetz index},
and $f$ has a global {\em Lefschetz number}. The classical Lefschetz
fixed point formula says that the sum of the Lefschetz indices is equal to the
Lefschetz number of $f$, up to a factor of $(-1)^n$
\cite{Lefschetz,Brown}.

In \cite[Lemma 3.7]{GoldbergMilnor}, Goldberg and Milnor give a
version of this theorem for weakly polynomial-like mappings
$f:\diskbar\to\C$. We prove a similar result for a class of maps
$f:\Delta\to\Cc$, where $\Delta\subset\Cc$ is a closed topological disk. By
extending the range of $f$ to $\Cc$, we allow poles and have to take more
boundary components into account than Goldberg and Milnor.

\begin{definition}[Lefschetz Map]
\label{Def_LefschetzMap} Let $\Delta\subset\Cc$ be a closed
topological disk with boundary curve $\partial \Delta$ and
$f:\Delta\to\Cc$ an orientation preserving open mapping with
isolated fixed points. We call $f$ a {\em Lefschetz map} if it
satisfies the following conditions: 
\begin{itemize}
\item
for every $z\in\Cc$, the full preimage $f^{-1}(\{z\})\subset\Delta$ is a finite set;
\item
 $f(\partial\Delta)$ is a simple closed curve so that
$f|_{\partial\Delta}:\partial\Delta\to f(\partial\Delta)$
is a covering map of finite degree;
\item
\(
    f(\partial\Delta)\cap \mathring{\Delta}=\emptyset
\)
\item
if $\xi\in\partial\Delta$ is a fixed point of $f$, then
$\xi$ has a neighborhood $U$ such that $f(\partial\Delta\cap U)\subset
\partial\Delta$, and $f$ is expanding on $\partial\Delta\cap U$.
\end{itemize}
\end{definition}
\begin{remark} The definition of ``expanding'' is with respect to the
local para\-me\-tri\-za\-tion of $\partial \Delta$ near $\xi$ so
that $f|_{\partial\Delta\cap U}$ is topologically conjugate to
$x\mapsto 2x$ in a neighborhood of $0$.

In this case, the map $f$ can be extended continuously to $U\cup\Delta$ so that
$f$ on $U\sm \mathring\Delta$ is topologically conjugate to $z\mapsto 2z$ on the
half disk $\{z\in\C\colon |z| < 1\mbox{ and } \Im(z)\ge 0\}$ (possibly after
shrinking $U$). Such an extension will be called the {\em simple extension outside
of $\Delta$} near $\xi$.
\end{remark}
\begin{definition}[Lefschetz Index]
\label{Def_Lefschetz}
Let $W\subset\C$ be a closed topological disk and $f\colon W\to \C$ be continuous
with an isolated fixed point at $\xi\in \mathring W$.  With $g(z)=f(z)-z$, we assign
to
$\xi$ its {\em Lefschetz index}
\[
    \iota(\xi,f) := \lim_{\eps\searrow 0}\; \frac{1}{2\pi i}\oint_{g(\partial B_{\eps}(\xi))} \frac{d\zeta}{\zeta}\;\;.
\]
This is the number of full turns that the vector $f(z)-z$ makes
when $z$ goes once around $\xi$ in a sufficiently small
neighborhood.

If $\xi\in\partial W\cap\C$ is an isolated boundary fixed point
which has a simple extension outside of $W$, then we define its
Lefschetz index as above for this simple extension.
\end{definition}
For an interior fixed point, it is easy to see that the limit
exists and is invariant under homotopies of $f$ that avoid
additional fixed points. Strictly speaking, the curve $g(\partial
B_\eps(\xi))$ need not be an admissible integration path (i.e.\ rectifiable), but
because of homotopy invariance, we may ignore this problem, and we will often do so
in what follows.

The Lefschetz index is clearly a local topological invariant; for boundary fixed
points, it does not depend on the details of the extension. Therefore, the index
is also defined if $\xi=\infty$, using local topological coordinates. Note that for
boundary points, the simple extension as defined above generates the least possible
Lefschetz index for all extensions of $f$ to a neighborhood of $\xi$.

If $f$ is holomorphic in a neighborhood of a fixed point
$\xi$, then $\iota(\xi,f)$ is the multiplicity of $\xi$ as a fixed
point.

\begin{definition}[Lefschetz Number]
\label{Def:LefschetzNumber}
Let $f\colon\Delta\to\Cc$ be a Lefschetz map, let $V$ be the component of $\Cc\sm f(\partial\Delta)$
containing $\mathring\Delta$ and let $\gamma_k$ denote the components of $\partial f^{-1}(V)$.
The {\em Lefschetz number} $L(f)$ of $f$ is then defined as
\[
L(f):=\sum_k \left|\deg\left(f|_{\gamma_k}\colon \gamma_k\to\gamma\right)\right| \,\,.
\]
\end{definition}

\begin{remark}
In this definition, the orientations of all $\gamma_k$ are irrelevant. The $\gamma_k$
are exactly the components of $f^{-1}(f(\partial\Delta))$, possibly with the exception of
$\partial\Delta$ itself: this latter curve is counted only if points in $\mathring\Delta$ near the
boundary of $\Delta$ are mapped into $V$ (in the other case, one can imagine
$\Cc\sm\ovl V$ as an omitted component of $f^{-1}(V)$, and consequently we also omit its boundary
curve). As a result, $L(f)\ge 0$, with equality iff $f^{-1}(V)\cap\Delta=\emptyset$, i.e., 
the sum over the $\gamma_k$ is empty.
\end{remark}

The Lefschetz number is invariant under topological
conjugacies. We may thus choose coordinates so that $\ovl V\subset\C$.

\begin{lemma}[Mapping Degree on Curve]
\label{Lem_WindingNumber} 
Let $\gamma\subset\C$ be a Jordan curve and $f\colon\gamma\to\C$ continuous
so that $f\colon\gamma\to f(\gamma)$ is a covering map and
 $\gamma\sm f(\gamma)$ is contained in the bounded component of $\C\sm
f(\gamma)$. If $\gamma$ contains no fixed points
of $f$, then the mapping degree of $f|_\gamma\colon\gamma\to f(\gamma)$ satisfies
\begin{equation}
	\deg(f|_\gamma\colon\gamma\to f(\gamma))=
    \frac{1}{2\pi i} \oint_{g(\gamma)}
    \frac{d\zeta}{\zeta}\;\;,
\label{Eq:WindingIntegral}
\end{equation}
where $g(z)=f(z)-z$, and $\gamma$ and $f(\gamma)$ inherit their orientations from $\C$.
\end{lemma}
\begin{proof} 
Let $\Delta$ be the bounded component of $\C\sm\gamma$ and 
$w_0\in\mathring{\Delta}$ any base point. Since
$\Delta$ is contractible to $w_0$ within $\mathring\Delta$,
$g|_{\partial\Delta}=(f-\id)|_{\partial\Delta}$ is homotopic to
$(f-w_0)|_{\partial\Delta}$ in $\C\sm\{0\}$.

The integral (\ref{Eq:WindingIntegral}) counts the number of full
turns of $f(z)-z$ as $z$ runs around $\gamma=\partial\Delta$. By homotopy invariance, this
is equal to the number of full turns $f(\gamma)$ makes
around $w_0$, and this equals the desired mapping degree of $f|_\gamma$.
\end{proof}

\begin{lemma}[Equality in $\C$]
\label{Lem:LocalIndex} Let $V\subset\C$ be a simply connected and
bounded domain with piecewise $\mathcal{C}^1$ boundary and let
$f\colon \ovl V\to f(\ovl V)\subset\C$ be a continuous map with
finitely many fixed points, none of which are on $\partial V$.
Then
\[
\sum_{f(\xi)=\xi} \iota(\xi,f) = \frac{1}{2\pi i}\oint_{g(\partial
V)}\frac{d\zeta}{\zeta}\;,
\]
where again $g=f-\id$.
\end{lemma}
\begin{proof} Break up $V$ into finitely many disjoint simply connected
open pieces $V_j$ with piecewise $\mathcal{C}^1$ boundaries so
that each $V_j$ either contains a single fixed point of $f$ or
$f(\ovl V_j)\cap \ovl V_j=\emptyset$, and each fixed point of $f$
is contained in some $V_j$. This can be done by first choosing
disjoint neighborhoods for all fixed points and then partitioning
their compact complement in $\ol{V}$ into pieces of diameter less
than $\theta$, where $\theta$ is chosen in such a way that
$|f-\id|>\theta$ in this complement. Set
\[
c_j:=\frac{1}{2\pi i}\oint_{g(\partial V_j)}\frac{d\zeta}{\zeta}
\;.
\]
Then
\[
\frac{1}{2\pi i}\oint_{g(\partial V)}\frac{d\zeta}{\zeta}
=
\sum_j \frac{1}{2\pi i}\oint_{g(\partial V_j)}\frac{d\zeta}{\zeta}
=
\sum_j c_j
\;.
\]
On the pieces with $f(\ovl V_j)\cap \ovl V_j=\emptyset$, we have $c_j=0$, and on a
piece $V_j$ with fixed point $\xi_j$, we have $\iota(\xi_j,f)=c_i$ by definition.
The claim follows.
\end{proof}

\begin{theorem}[Fixed Point Count]
\label{Thm:Lefschetz} 
Let $f:\Delta\to\Cc$ be a Lefschetz map with
Lefschetz number $L(f)\in\N$. Then
\[
    L(f) = \sum_{f(\xi)=\xi} \iota(\xi,f)\;\;.
\]
\end{theorem}
\begin{proof} Let $V$ be the component of $\Cc\sm f(\partial\Delta)$
containing $\mathring\Delta$, and choose coordinates of $\Cc$ such
that $\ovl V$ is bounded.

Suppose first that $f$ has no fixed points on $\partial\Delta$.
Let $\{U_i\}$ be the collection of components of $f^{-1}(\Cc\sm
\ovl V)$ and let $\{V_j\}$ be the collection of components of
$f^{-1}(V)$. Since $f$ is open, each $U_i$ maps onto $\Cc\sm \ovl
V$ and each $V_j$ maps onto $V$ as a proper map. It follows that
there are only finitely many $U_i$ and $V_i$, and they satisfy
$f(\partial U_i)=f(\partial V_i)=f(\partial\Delta)$ for each $U_i$ and each $V_j$.
Note that every fixed point of $f$ must be in some $V_j$.

Subdivide the $V_j$ into finitely many simply connected pieces so that no
fixed points of $f$ are on the boundaries; call these subdivided domains
$V'_{j'}$. The orientation of $\C$ induces a boundary orientation on the
$V'_{j'}$.

Set again $g:=f-\id$. Then, applying Lemma~\ref{Lem:LocalIndex} to
$V'_{j'}\subset\Delta\subset\ovl{V}\subset\C$ yields
\[
\sum_{f(\xi)=\xi}\iota(\xi,f)
=
\sum_{j'}
\frac{1}{2\pi i} \oint_{g(\partial V'_{j'})}\frac{d\zeta}{\zeta}
=
\sum_{j}
\frac{1}{2\pi i} \oint_{g(\partial V_j)}\frac{d\zeta}{\zeta}
\;.
\]
Let $\gamma:=\partial V$ and  $\gamma_k \subset\Delta$ be the components of 
$\partial V_k$ (these are exactly the curves from Definition~\ref{Def:LefschetzNumber}), with the orientation they inherit from $\partial V_k$.
Then
\[
\sum_{j}
\frac{1}{2\pi i} \oint_{g(\partial V_j)}\frac{d\zeta}{\zeta}
=
\sum_k
\frac{1}{2\pi i} \oint_{g(\gamma_k)}\frac{d\zeta}{\zeta}
\,\,.
\]
The curves $\gamma_k$ come in two different kinds: those which surround the
component $V_j$ of which they form part of the boundary, and those which are surrounded
by their component $V_j$. The first kind has the same orientation as the orientation
it would inherit as a simple closed curve in $\C$, and the second kind has the opposite
orientation. On the other hand, it is easy to check that $f|_{\gamma_k}\colon\gamma_k\to\gamma$
has positive mapping degree (with respect to this orientation of $\gamma_k$, and the 
standard orientation in $\C$ for $\gamma$) exactly for curves $\gamma_k$ of the first kind.
As a result, Lemma~\ref{Lem_WindingNumber} implies in both cases
\[
\frac{1}{2\pi i} \oint_{g(\gamma_k)}\frac{d\zeta}{\zeta}
=
\left|\deg\left(f|_{\gamma_k}\colon \gamma_k\to\gamma\right)\right|
\,\,.
\]
This implies the claim if $f$ has no fixed points on $\partial\Delta$.

If $f$ has boundary fixed points, we employ a simple extension
outside of $\Delta$ in a small neighborhood of each such fixed
point. In order for the extended map to be a Lefschetz map, the
preimages need to be extended as well. If the extended
neighborhoods are sufficiently small, this does not change the
Lefschetz number of $f$.
\end{proof}

As an immediate corollary of this theorem, we observe that Newton maps do not have fixed Herman rings. Note that Taixes has announced a more general result: he uses quasiconformal surgery to rule out any periodic cycles of Herman rings for Newton maps. 

\begin{corollary}[No Fixed Herman Rings]
\label{Cor_Herman}
Newton maps of entire functions have no fixed Herman rings.
\end{corollary}
\begin{proof}
By \cite{Shishikura}, we may assume that $N:\C\to\Cc$ is a transcendental meromorphic Newton map. 
Suppose it has a fixed Herman ring, i.e.\ an invariant Fatou component $H$ such that $N|_H$ is conjugate to an irrational rotation of an annulus of finite modulus. Then, $H$ contains an invariant and essential simple closed curve $\gamma$. Clearly, ${\rm deg}(N:\gamma\to\gamma)=+1$. Let $\Delta$ be the bounded component of $\C\setminus\gamma$. Then, $N|_{\ol{\Delta}}$ is a Lefschetz map and by Theorem \ref{Thm:Lefschetz}, $\ol{\Delta}$ contains a fixed point. This is a contradiction, because all fixed points of $N$ have an unbounded immediate basin (Theorem \ref{Thm_ImmediateBasinSC}).
\end{proof}

\section{Between Accesses of an Immediate Basin}
\label{Sec_Main}

In this section, we state and prove our main result. Let
$f:\C\to\C$ be an entire function and $N_f$ its Newton map. Let
$\xi\in\C$ be a fixed point of $N_f$ and $U$ its immediate basin.
Suppose that $U$ has two distinct invariant accesses, represented
by $N_f$-invariant curves $\Gamma_1$ and $\Gamma_2$. Consider an
unbounded component $\tilde{V}$ of $\C\sm (\Gamma_1\cup\Gamma_2)$.
We keep this notation for the entire section.

\begin{theorem}[Main Theorem]
\label{Thm_NecessaryCondition}
If no point in $\Cc$ has infinitely many preimages within $\tilde V$, then
the set $V:=\tilde{V}\sm U$ contains an immediate basin or a virtual
immediate basin of $N_f$.
\end{theorem}
Note that we do not assume that $V$ is connected.
\begin{corollary}[Polynomial Case]
\label{Cor_PolynomialCase}
    If $N_f$ is the Newton map of a polynomial $f$, then each
    component of $\C\sm U$
    contains the immediate basin of another root of $f$.
\end{corollary}
\begin{proof}[Proof of Corollary~\ref{Cor_PolynomialCase}] 
If $f$ is a polynomial, $N_f$ is a rational map. It has finite mapping degree
and there exists $R>0$ such that all components of $U\sm B_R(0)$
contain exactly one invariant access. Furthermore, all accesses
are invariant \cite[Proposition 6]{HSS}. Since $\infty$ is a
repelling fixed point of $N_f$, there are no virtual immediate
basins.
\end{proof}
The rest of this section will be devoted to the proof of Theorem
\ref{Thm_NecessaryCondition}. This proof will be based on the
fixed point formula in Theorem \ref{Thm:Lefschetz}. In order to be
able to use it in our setting, we will need some preliminary
statements.
\begin{proposition}[Pole on Boundary]
\label{Pole_on_Boundary}
    If no point $z\in U$ has infinitely many preimages within
    $\tilde{V}\cap U$, then $\partial V=\partial U\cap \tilde{V}$
    contains at least one pole of $N_f$.

    In particular, if
    $\partial V$ is connected or $N_f|_U$ has finite degree,
    then $\partial V$ contains a pole of $N_f$.
    Every pole on $\partial V$ is arcwise accessible from within $U$.
\end{proposition}
\begin{proof} Let $\phi:\disk\to U$ be a Riemann map for the immediate
basin $U$ with $\phi(0)=\xi$. It conjugates the dynamics of $N_f$
on $U$ to the induced map $g=\phi^{-1}\circ
N_f\circ\phi:\disk\to\disk$. By Lemma \ref{Lem_RiemannConnected},
the Carath\'eodory extension $\ol{\phi}^{-1}$ maps $\ol{\partial
V}\subset\partial U\cup\{\infty\}$ to a closed interval $I\subset
\partial\disk$ that is bounded by
the landing points $\zeta_1$ and $\zeta_2$ of
$\phi^{-1}(\Gamma_1)$ and $\phi^{-1}(\Gamma_2)$. By assumption,
there is an open neighborhood of $\mathring{I}$ in $\disk$ which contains
only finitely many $g$-preimages of every $z\in\disk$. By Proposition
\ref{PropInvAccesshasFiniteDegree}, there is a neighborhood $W'$ of $I$ in
$\disk$ with the same property. Consider a sequence $(z_n)\subset\disk$ whose
accumulation set is in $W'\cap\partial\disk$. By Lemma \ref{No_Asymptotic_Path}, all
limit points of $(g(z_n))$ are in $\partial\disk$. Hence there is
a neighborhood $W$ of $I$ in $\C$ such that we can extend $g$ by
Schwarz reflection to a holomorphic map $\tilde{g}:W\to\C$ that
coincides with $g$ on $W\cap\disk$. The endpoints $\zeta_1$ and
$\zeta_2$ of $I$ are fixed under this map, because each is the
landing point of an invariant curve. They are repelling, because
otherwise they would attract points from within $\disk$.

Clearly, $\tilde{g}(I)\subset\partial\disk$. If $\tilde{g}(I)=I$,
then $\tilde{g}$ has to have an additional fixed point on
$\mathring{I}$ which is necessarily parabolic and thus attracts
points in $\disk$. This is a contradiction because all points in
$\disk$ converge to $0$ under iteration of $g$.

If $I$ contained a critical point $c$ of $\tilde{g}$, points in
$\disk$ arbitrarily close to $c$ would be mapped out of $\disk$ by
$\tilde{g}$, again a contradiction. Hence
$\tilde{g}:I\to\partial\disk$ is surjective and there are points
$z_1, z_2\in\mathring{I}$ such that
$\tilde{g}(z_1)=\zeta_1,\,\tilde{g}(z_2)=\zeta_2$.

For $i=1,2$, let $\beta_i:[0;1)\to\disk$ be the radial line from
$0$ to $z_i$. Then, $\phi(\beta_i)$ accumulates at a continuum
$X_i\subset \partial V$ while
$N_f(\phi(\beta_i))=\phi(g(\beta_i))$ lands at $\infty$ in the
access of $\Gamma_i$. By continuity, $N_f(X_i)=\{\infty\}$; the
identity theorem shows that $X_i=\{p_i\}$ is a pole and
$\phi(\beta_i)$ lands at $p_i$.
\end{proof}
We use the following general lemma to show that $N_f|_{\tilde{V}}$
can be continuously extended to $\infty$.
\begin{lemma}[Extension Lemma]
\label{contExtension}
    Let $h:\C\to\Cc$ be a meromorphic function and $G\subset \C$ an unbounded
    domain. Suppose that $\partial G$ can be parametrized
    by two asymptotic paths of the asymptotic value $\infty$ and
    that no point has infinitely many preimages within $G$.
    Then, $h|_{\ol{G}}$ can be continuously extended to $\infty$.
\end{lemma}
\begin{proof} Since $h(\partial G)$ is unbounded, the only possible
continuous extension is to set $h(\infty)=\infty$.

If $h$ cannot be continuously extended to $\infty$, there exists a
sequence $z_n\to\infty$ in $G$ such that $h(z_n)\to p\in\C$. Let
$S>|p|$ and pick $R>0$ such that $|h(z)|\geq S$ for all
$z\in\partial G$ with $|z| \geq R$, and $p\notin h(\partial
B_R(0))$. We may suppose that all $|z_n|>R$. Then we can choose a
closed neighborhood $B\subset B_S(0)$ of $p$ whose boundary is a
simple closed curve that contains no critical values or direct
singularities and so that $B$ is disjoint from $h(\partial
B_R(0))$. Now let $W_n$ be the component of $h^{-1}(B)$ that
contains $z_n$. Then $W_n\subset \C\sm \ovl{B_R(0)}$. Since
$z_n\in G$, it follows that all $W_n\subset G\sm \ovl{ B_R(0)}$.

If all $W_n$ are bounded, each can contain only finitely many
$z_k$ and there must be infinitely many such components. Since
bounded $W_n$ map onto $B$, this would contradict the finiteness
assumption. Hence there is an unbounded preimage component $W_0$.
By Lemma \ref{Lem_AsymptoticValueOnBoundary}, $G\sm \ol{B_R(0)}$
then contains an asymptotic path of an indirect singularity on
$\partial B$, which also contradicts the finiteness assumption.
\end{proof}
In the next proposition, we show that $N_f|_{\tilde{V}}$ is
injective near $\infty$. For the proof, we use an extremal length
argument in the half-strip
\[
Y:=[0,\infty)\times [0,1]\;\;,
\]
in which we measure the modulus of a quadrilateral by curves
connecting the left boundary arc to the right. For $x\in\R$,
define
\[
\H_{x}:=\{z\in\C\,:\,\Re(z)\geq x\}\;\;.
\]
First, we prove a technical lemma.
\begin{lemma}[Bound on Modulus]
\label{Lem_BoundedModulus} Let $0<t\le s$, let $\beta\subset Y$ an
injective curve from $(t,1)$ to $(s,0)$ and let $Q$ the bounded
component of $Y\sm\beta$. Let $(0,0)$, $(s,0)$, $(t,1)$ and
$(0,1)$ be the vertices of the quadrilateral $Q$. Then, $\mod(Q)
\leq t+1$.
\end{lemma}
\begin{proof} Let $R\subset Y$ be the rectangle with vertices $(0,0)$,
$(0,1)$, $(t+1,1)$, $(t+1,0)$. Its area and modulus are both equal
to $t+1$. In particular, $\area(Q\cap R)\leq\area(R)=t+1$. Using
the admissible density $\rho(x)=\frac{1}{\sqrt{\area(Q\cap
R)}}\cdot\chi_{Q\cap R}(x)$, we get the estimate
\[
    \frac{1}{\mod(Q)}\geq\frac{1}{\area(Q\cap R)} \geq
    \frac{1}{\area(R)}=\frac{1}{t+1}\;\;,
\]
because $\int_\gamma \rho\, d\gamma\geq \frac{1}{\sqrt{\area(Q\cap
R)}}$ for this density and all rectifiable curves $\gamma$ that
connect the upper to the lower boundaries.
\end{proof}
\begin{proposition}[Invariance Near $\infty$]
\label{Prop_InvarianceNearInfinity} Suppose that every $z\in\Cc$
has only finitely many $N_f$-preimages in $\tilde V$. Then there
exists $R_0>0$ such that for all $R>R_0$, the map $N_f$ is
injective on $\tilde{V}\sm B_R(0)$. Moreover, there exists $S>0$
with the property that
\[
N_f(\tilde{V}\sm B_R(0))\sm B_S(0)=\tilde{V}\sm
B_S(0) \,\,.
\]
\end{proposition}
\begin{proof}
Choose $R_0>\max\{|z|:z\in N_f^{-1}(\infty)\cap
\ol{\tilde{V}}\}$. It follows from the open mapping principle and
invariance of $\partial\tilde{V}$ that there exists $S_0>0$ such
that $\partial N_f(\tilde{V}\sm B_{R_0}(0))\sm
B_{S_0}(0) \subset\partial\tilde{V}$.
Since there are points $z\in\tilde V$ with arbitrarily large $|z|$
such that $N_f(z)\in\tilde{V}$, it follows
that either $N_f(\tilde{V}\sm B_R(0))\subset \tilde{V}$
or $N_f(\tilde{V}\sm B_R(0))$ contains a punctured neighborhood
of $\infty$ within $\Cc$.

In the first case, the claims follow easily.
By way of contradiction, we may thus assume that we are in the second case.

We consider the situation in logarithmic coordinates: with an
arbitrary but fixed choice of branch, let $C\subset\H_{\log(R)}$
be the unique unbounded component of $\log(\ol{\tilde{V}}\sm
B_R(0))$. This is a closed set whose boundary consists of two
analytic curves $\gamma_1$ and $\gamma_2$ and a subset of the
vertical line at real part $\log(R)$. Define a holomorphic map
$g:C\to\C$ by $g(z)=\log(N_f(e^z))$, choosing the branch such that
$\gamma_1\subset g(\gamma_1)$. This is possible because
$\Gamma_1=e^{\gamma_1}$ is $N_f$-invariant. Since $e^{\gamma_2}$
is also $N_f$-invariant, there exists $k\in\Z$ such that with
$\gamma_4:=g(\gamma_2)$, $\gamma_4=\gamma_2+2\pi i k$. Define also
$\gamma_3:=\gamma_1+2\pi i k$. Since $N_f(\tilde{V}\sm B_R(0))$
contains a neighborhood of $\infty$, we get $k\neq 0$. See Figure
\ref{Figure_Injectivity} for an illustration of the notations and
note that $\gamma_1,\gamma_2,\gamma_3$ and $\gamma_4$ are pairwise
disjoint. These four curves have a natural vertical order induced
by the observation that each curve separates sufficiently far
right half planes into two unbounded components. To fix ideas,
suppose that $\gamma_2$ is below $\gamma_1$. Then $\gamma_4$ is
below $\gamma_3$. The construction implies that $\gamma_4$ is
below $\gamma_1$, and no curve is between $\gamma_1$ and
$\gamma_2$. Then, the vertical order is $\gamma_1$, $\gamma_2$,
$\gamma_3$, $\gamma_4$.
\begin{figure}[ht]
    \begin{center}
        \setlength{\unitlength}{1cm}
        \begin{picture}(12,5.3)
        \epsfig{file=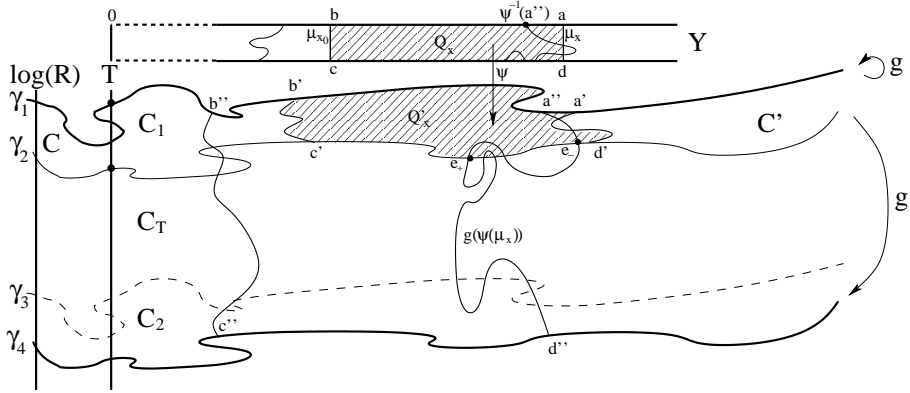, width=12cm}
        \end{picture}
        \caption{Illustrating the notations for Proposition \ref{Prop_InvarianceNearInfinity}.}
        \label{Figure_Injectivity}
    \end{center}
\end{figure}

By Lemma \ref{contExtension}, $g$ is continuous at $+\infty$. By
the open mapping principle, there exists $T>\log(R)$ such that
$\partial g(C)\cap \H_T \subset \gamma_1\cup\gamma_4$. Let $C_1$
be the unique unbounded component of $C\cap \ol{\H}_T$,
$C_2=C_1+2\pi i k$ and $C_T := g(C)\cap \ol{\H}_T$. Note that
$C_1\cup C_2\subset C_T$. We may choose $T$ in such a way that
$\partial \H_T$ does not contain any critical values of $g$.
Define $C'=g^{-1}(C_T)\subset C$. Then, $g:C'\to C_T$ is a proper
map and therefore has well-defined degree. Since $g$ is injective
on $\partial C'$ and has no pole, this degree is one and $g$ is
univalent.

The idea of the proof is as follows: the curves $\gamma_2$ and
$\gamma_3$ subdivide $C_T$ into three parts which are unbounded to
the right. With an appropriate bound to the right, we obtain a
large bounded quadrilateral consisting of three
sub-quadrilaterals. Two of these sub-quadrilaterals, the upper and
the lower ones, have moduli comparable to the modulus of the
entire quadrilateral. This is a contradiction to the Gr\"otzsch
inequality if the right boundaries are sufficiently far out.

Define a homeomorphism $\psi:Y\to C_1$ that is biholomorphic on
the interior and normalized so that it preserves the boundary
vertex $\infty$ and the other two boundary vertices. We denote by
$\mu_x\subset Y$ the vertical line segment at real part $x$. There
exists an $x_0$ such that for $x\geq x_0$, $\psi(\mu_{x})\subset
C'$. For $x>x_0$, we denote by $Q_x$ the rectangle in $Y$ that is
bounded by $\mu_{x_0}$ and $\mu_x$. With vertices $a=(x,1)$,
$b=(x_0,1)$, $c=(x_0,0)$ and $d=(x,0)$, its modulus is equal to
$x-x_0$. We denote the vertices of its image $Q'_x:=\psi(Q_x)$ by
$a',b',c',d'$, respectively. Let $a''=g(a')$ and $d''=g(d')$.
Since $g$ and $\psi$ are univalent, $\mod(g(Q_x'))=x-x_0$.

The curve $g(\psi(\mu_x))$ is a boundary curve of $g(Q'_x)$; it
connects $a''$ and $d''$ within $C_T$. Let $e_-$ be the
intersection point $g(\psi(\mu_x))\cap\gamma_2$ closest to $a''$
along $g(\phi(\mu_x))$, and let $e_+$ be the intersection point
furthest to the left along $\gamma_2$. Let $C'_1$ be the bounded
subdomain of $C_1$ bounded by $g(\psi(\mu_x))$ between $a''$ and
$e_-$, viewed as a quadrilateral with vertices $a''$ and $e_-$ and
two more vertices on $\partial\H_T$. Similarly, let $C''''_1$ be
the bounded subdomain of $\C$ bounded by $\partial\H_T$,
$\gamma_1$, the part of $\gamma_2$ to the left of $e_+$, and the
part of $g(\phi(\mu_x))$ between $a''$ and $e_+$, with right
vertices $a''$ and $e_+$. Finally, let $C''_1:=C'_1\cup C''''_1$
with right vertices $a''$ and $e_-$, and let $C'''_1:=C''_1$ but
with right vertices $a''$ and $e_+$ (instead of $a''$ and $e_-$).

If $g(\psi(\mu_x))$ intersects $\gamma_2$ only once, then
$e_-=e_+$ and $C'_1=C''_1=C'''_1=C''''_1$. In general, the three
domains $C'_1,C''_1,C''''_1$ may be different. However, we have
$\mod(C'_1)\ge\mod(C''_1)\ge\mod(C'''_1)\ge\mod(C''''_1)$: the
first inequality holds because $C'_1\subset C''_1$, the second
describes identical domains but with one boundary vertex moved,
and the third follows again from the inclusion $C''''_1\subset
C'''_1$, but this time the domain is extended on the ``right''
side of the domain, rather than on the ``lower'' side because the
boundary vertex has moved.

Pulling back under $\psi$, we find that $\Re(\psi^{-1}(a''))\le
\Re(a)=x$, because the map $\psi^{-1}\circ g\circ\psi$ repels
points away from $\infty$. By Lemma~\ref{Lem_BoundedModulus}, it
follows that $\mod(C''''_1)\le\mod(C'_1)\le x+1$.

Similar considerations on the left end of $C_1$, as well as for $C_2$, allow to
subdivide $g(Q'_x)$ by a single curve segment of $\gamma_2$ and $\gamma_3$ into
three sub-quadrilaterals, two of which have modulus at most $x+1$. But the
Gr\"otzsch inequality implies that
\[
    \frac{1}{x-x_0}=\frac{1}{\mod(g(Q_x'))}\geq \frac{1}{x+1}+\frac{1}{x+1}\;,
\]
hence $x\le 2x_0+1$ which is a contradiction for large $x$.
\end{proof}
\begin{proof}[Proof of Theorem~\ref{Thm_NecessaryCondition}] 
In order to use
Theorem \ref{Thm:Lefschetz}, we  construct an injective curve that
surrounds an unbounded domain in $\tilde{V}$ such that the image
of the curve does not intersect this domain. Consider a Riemann
map $\phi:\disk\to U$ with $\phi(0)=\xi$ and the induced dynamics
$g=\phi^{-1}\circ N_f\circ\phi$ on $\disk$. By the Remark after
Proposition \ref{PropInvAccesshasFiniteDegree}, the curves
$\phi^{-1}(\Gamma_1)$ and $\phi^{-1}(\Gamma_2)$ land at points
$\zeta_1,\zeta_2\in\partial\disk$, and $g$ extends to a
neighborhood of $\zeta_1$ and $\zeta_2$ so that $\zeta_1$ and
$\zeta_2$ become repelling fixed points. These fixed points have
linearizing neighborhoods in which the curves
$\phi^{-1}(\Gamma_1)$, respectively $\phi^{-1}(\Gamma_2)$, are
straight lines in linearizing coordinates. If $0<r<1$ is large
enough, these two curves intersect the circle at radius $r$ only
once and we can join them by a circle segment at radius $r$ to an
injective curve $\Gamma'\subset\disk$ in such a way that
$\Gamma:=\phi(\Gamma')$ separates $V$ from $\xi$. Let $W$ be the
closure in $\Cc$ of the connected component of $\C\sm\Gamma$ that
contains $V$ (Figure \ref{Figure_LastProof}). Note that no
component of $N_f^{-1}(\Gamma)$ that intersects $\mathring{W}$ can
leave $W$: in $\disk$, any such component would have to intersect
$\Gamma'$. But by the Schwarz Lemma, $g^{-1}(\Gamma')$ has greater
absolute value than $r$ everywhere and $\Gamma'$ has only one
$g$-preimage within the linearizing neighborhood of $\zeta_1$;
this preimage is contained in $\Gamma'$. The same is true at
$\zeta_2$.
\begin{figure}[ht]
    \begin{center}
        \setlength{\unitlength}{1cm}
        \begin{picture}(8,6)
            \epsfig{file=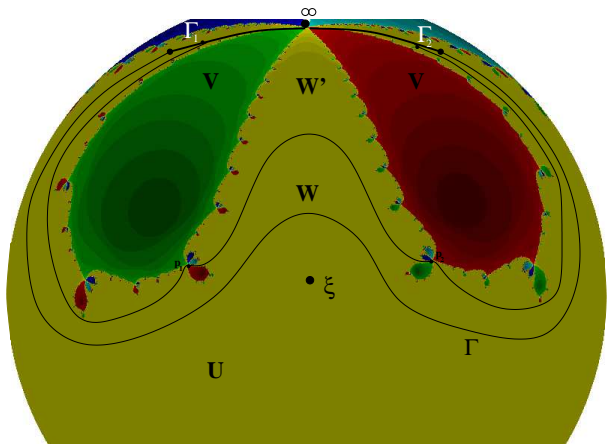, width=8cm}
        \end{picture}
        \caption{A schematic illustration of some notations in the proof of
Theorem~\ref{Thm_NecessaryCondition}.}
        \label{Figure_LastProof}
    \end{center}
\end{figure}

By Proposition \ref{Prop_InvarianceNearInfinity}, $W$ contains an
unbounded preimage component $W'$ of itself such that the boundary
$\partial W'$ is contained in $\Gamma_1\cup\Gamma_2$ outside a
sufficiently large disk. Make $W'$ simply connected by filling in
all bounded complementary components.

We claim that $\partial W'$ contains at least one finite pole on
$\partial U$: if it did not, then $\partial W'\subset U$ and
$N_f|_{\partial W'} : \partial W' \to \Gamma$ would be injective
for all choices of $r$ above. In the limit for $r\to 1$, this
would imply that $N_f|_{\partial V}$ was injective, contradicting
Proposition~\ref{Pole_on_Boundary}.

Therefore, $\partial W'$ maps onto $\Gamma$ with covering degree at
least $2$. If $\infty$ is not an isolated fixed point in $W'$, we
are done. Otherwise it is easy to see that $N_f|_{W'}$ is a
Lefschetz map: there is a single boundary fixed point $\infty$;
the conditions on this boundary fixed point are satisfied because
$\Gamma_1,\Gamma_2\subset U$, where the dynamics is expanding away
from $\infty$. Now Theorem~\ref{Thm:Lefschetz} implies that $W'$
contains fixed points of combined Lefschetz indices at least $2$, 
because $\partial W'$ contains a pole. If $W'$ contains a finite
fixed point, we are done. If not, it follows that the fixed point
at $\infty$ has Lefschetz index at least $2$. Consider a Riemann
map $\psi: W' \to \H^+$ that uniformizes $W'$ to the upper half
plane and maps $\infty$ to $0$; this map preserves the Lefschetz
index. By Proposition~\ref{Prop_InvarianceNearInfinity}, the map
$g=\psi\circ N_f\circ\psi^{-1}$ is defined in a relative
neighborhood of $0$ in $\H^+$. If a sequence converges to $\R$ in
this neighborhood, then so will the image of this sequence. Hence
we can extend $g$ to a neighborhood of $0$ in $\C$ by reflection.
This extension does not reduce the Lefschetz index of $0$: for a
boundary fixed point, the index is defined by extending $g$ to the
lower half-plane in the way which generates the least possible
fixed point index (compare Definition~\ref{Def_Lefschetz}).
Reflection however may increase the Lefschetz index. Therefore,
$0$ is a parabolic (since multiple) fixed point of the extended
map, and it is easily seen that $\partial\H^+$ is in the repelling
direction. By the Fatou flower theorem~\cite[Theorem
10.5]{Milnor}, $0$ has an attracting petal in $\H^+$ that induces
a virtual immediate basin inside $V$.
\end{proof}
\section*{Acknowledgements}
We would like to thank Walter Bergweiler, Xavier Buff, Matthias
G\"orner, Alexandra Kaffl, Sebastian Mayer and Lasse Rempe for
their helpful comments in many discussions and their support. We
also thank the Institut Henri Poincar\'e, Universit\'e Paris VI,
where much of the paper evolved.

\end{document}